\documentclass[11pt]{article}

\usepackage{amssymb,amsmath,amscd,amsthm,epsfig}
\usepackage{amsfonts}
\usepackage{latexsym}
\newcommand{\bbz}{{\mathbf Z}}%{\mbox{\boldmath $Z$}}
\newcommand{\bbq}{{\mathbf Q}}
\newcommand{\bbf}{\mbox{\boldmath $F$}}

\newcommand{\Des}{{\rm Des}}

\newcommand{\maj}{{\rm maj}}

\newcommand{\fmaj}{{\rm fmaj}}
\newcommand{\Fmaj}{{\rm Fmaj}}

\newcommand{\bc}{{\bar{c}}}

\newcommand{\tk}{{\tilde{k}}}

\newcommand{\Hilb}{{\rm{Hilb}}}
\newcommand{\Poin}{{\rm{Poin}}}
\newcommand{\then}{{\Longrightarrow}}

\newtheorem{thm}{Theorem}[section]
\newtheorem{pro}[thm]{Proposition}
\newtheorem{lem}[thm]{Lemma}

\newtheorem{cor}[thm]{Corollary}
\newtheorem{fac}[thm]{Fact}
\newtheorem{obs}[thm]{Observation}

\newtheorem{df}[thm]{Definition}
\newtheorem{rem}[thm]{Remark}

\title%[Complex Reflection Groups]
{Major Indices and %Ordered Generating Systems
Perfect Bases\\ for Complex Reflection Groups}

\author{%
Robert Shwartz%
\thanks{Department of Mathematics, Bar-Ilan University,
Ramat-Gan 52900, Israel. Email: {\tt shwartr1@math.biu.ac.il}}
\and Ron M. Adin%
\thanks{Department of Mathematics, Bar-Ilan University,
Ramat-Gan 52900, Israel. Email: {\tt radin@math.biu.ac.il}}
\and Yuval Roichman%
\thanks{Department of Mathematics, Bar-Ilan University,
Ramat-Gan 52900, Israel. Email: {\tt yuvalr@math.biu.ac.il}}
\thanks{Research of all authors was supported in part by
the Israel Science Foundation, founded by the Israel Academy of
Sciences and Humanities.}}

\date{August 6, 2007}

\begin{document}

\maketitle

\begin{abstract}
%A classical result of MacMahon shows that the length function and
%the major index are equi-distributed over the symmetric group.
%A long-standing open problem is to extend the notion of {\it major index}
%and MacMahon's result to other groups. A partial solution was given
%in~\cite{AR} and~\cite{BC1}, where this result was extended to
%classical Weyl groups.
%%A group has an ordered generating system if it may be decomposed
%%into a product of "almost disjoint" cyclic subgroups (namely
%%subgroups whose intersection is the identity element).
%In this paper it is proved
It is shown that, under mild conditions, a complex reflection group
$G(r,p,n)$ may be decomposed into a set-wise direct product of
cyclic subgroups.
% a product of "almost disjoint" cyclic subgroups
%(namely subgroups whose intersection is the identity element).
This property is then used to extend the notion of major index and
a corresponding  Hilbert series identity to these and other
closely related groups.
\end{abstract}

\section{Introduction}

\subsection{The major index}

Let $S_n$ be the symmetric group on $n$ letters.
$S_n$ is a Coxeter group with respect to the Coxeter generating set
$S=\{ s_i\,|\,1\le i<n\}$, where $s_i:=(i,i+1)$ ($1\le i<n$) are
the adjacent transpositions.
Let $\ell(\pi)$ be the {\em length} of $\pi\in S_n$ with respect to $S$,
let $$\Des(\pi):=\{1\leq i< n |\ \ell(\pi s_i)<\ell(\pi)\}
$$
be the {\em descent set} of $\pi$, and let
$$
\maj(\pi):=\sum_{i\in \Des(\pi)}i
$$
be the {\em major index} of $\pi$. It is well known that
$$
\ell(\pi) = \#\{i<j\ |\ \pi(i)>\pi(j)\},
$$
the number of inversions in $\pi$, and that
$$
\Des(\pi)=\{1\leq i\le n-1\,|\,\pi(i)>\pi(i+1)\}.
$$
The major index is involved in many classical %combinatorial
identities on the symmetric group;
%such as MacMahon's~\cite{MM}, Foata-Sch\"utzenbenger's~\cite{FS}
%and others,
see, for example, \cite{MM, FS, GaG, BW}.
The search for an extended major index and corresponding identities
on other groups, initiated by Foata in the early nineties, turned out
to be successful for the classical Weyl groups and some wreath products.
In particular, the Hilbert series of the coinvariant algebra of
the symmetric group $S_n$ and of the wreath products $\bbz_r\wr S_n$
may be expressed as generating functions for the {\em flag major index}
on these groups~\cite{AR, BC1, ABR}.
A generalization of this result to complex reflection groups,
involving the notion of {\em basis} for a group, is suggested in this
paper. This generalization extends previous results of~\cite{AR}.

\subsection{Bases}%{Ordered Generating Systems}

The concept of basis for a group~\cite{S1, Mu} extends the classical
Fundamental Theorem for Finitely Generated Abelian Groups to the
non-abelian case.

\begin{df}%\label{d.ogs*}
Let $G$ be a finite group. A sequence ${\bf a}=(a_1,\ldots, a_n)$
of elements of $G$ is called a {\em basis} (or a {\em starred ordered
generating system, OGS*}) for $G$ if there exist positive integers
$m_1,\ldots,m_n$ such that every element $g\in G$ has a unique
presentation in the form
$$
g=a_1^{k_1}a_2^{k_2}\cdots a_n^{k_n},
$$
with $0\le k_i<m_i$ for every $1\le i\le n$.

If $m_i=o(a_i)$ (the order of the element $a_i$) for every $1\le i\le n$ then
${\bf a}$ is a {\em perfect basis} (or an {\em ordered generating system, OGS})
for $G$.
\end{df}

%Denote the order of an element $g\in G$ by $o(g)$.
%
%\begin{df}\label{ogs}%\begin{itemize}
%%\item[1.] The sequence  $A=(a_1,\ldots, a_n)$ is called a {\rm
%%starred Ordered Generating System (OGS$^*$)} for a group $G$, if
%%there exist positive integers $m_1,\dots,m_n$, such that every
%%element $g\in G$ has a unique presentation:
%%$$g=a_1^{k_1}a_2^{k_2}\cdots a_n^{k_n},$$ where $0\leq k_i<m_i$, for every
%%$1\leq i\leq n$. \item[2.]
%The sequence  ${\bf a}=(a_1,\ldots, a_n)$ is called a {\rm perfect
%basis} (or an {\rm Ordered Generating System (OGS)}) for a group
%$G$, if every element $g\in G$ has a unique presentation:
%\begin{equation}\label{e.unique-presentation}
%g=a_1^{k_1}a_2^{k_2}\cdots a_n^{k_n}, \end{equation} where $0\leq
%k_i<o(a_i)$, for every $1\leq i\leq n$.
%%\end{itemize}
%%\end{df}
%\end{df}

A finite group $G$ has a perfect basis if and only if $G$ has a
decomposition into a set-wise direct product of cyclic subgroups.
Namely, a group $G$ has a perfect basis if and only if there exist
subgroups $C_1,\ldots,C_n$ of $G$ such that
$$
C_i \hbox{\ is cyclic } \qquad (\forall i), \leqno(i)
$$
$$
G=C_1\cdots C_n, \hbox{\ and} \leqno(ii)
$$
$$
C_i \cap \left(C_1\cdots \hat C_i \cdots C_n\right) = \{1\} \qquad
(\forall i), \leqno(iii)
$$

%An ordered generating set which gives a unique expression of type
%(\ref{e.unique-presentation}) but does not give a decomposition into
%a set-wise direct product of the cyclic subgroups is called an
%OGS$^*$, see~\cite{S1}.

\noindent{\bf Examples:}
\begin{enumerate}
\item %1)
$pq$-groups ($p,q$ distinct primes) have a perfect basis~\cite{S1}.

\item %2)
The group of quaternions $Q_8$ has a basis, but
{\bf not} a perfect basis~\cite{S1}.
\end{enumerate}
%An example of a finitely generated group that has {\bf no} basis
%is given~\cite{S1}.

%\bigskip

%\begin{prb} Characterize finite groups with OGS.
%\end{prb}

%\bigskip

\bigskip

The major index of a permutation has an algebraic interpretation in
terms of a perfect basis. The following observation is a
reformulation of~\cite[Claim 2.1]{AR}.

\begin{obs}\label{t.A-fmaj}%\cite[Claim 2.1]{AR}
Let $s_i := (i,i+1)\in S_n$ $(1\le i< n)$ and
$$
t_i := s_i s_{i-1}\cdots s_1 \qquad (1\le i< n).
$$
%be $n-1$ elements in $S_n$.
Then $(t_{n-1},t_{n-2},\ldots,t_1)$ is a perfect basis for $S_n$;
namely, every permutation $\pi\in S_n$ has a unique presentation
$$
\pi={t_{n-1}}^{k_{n-1}}\cdots {t_1}^{k_1},
$$
where  $0\le k_i < o(t_i)=i+1$ $(1\le i< n)$. In this notation,
$$
\maj(\pi) = \sum_{i=1}^{n-1}k_i.
$$
%For every permutation $\pi\in S_n$
%$$maj(\pi)=\sum_{i=1}^{n-1}k_i,$$ where $\pi$ is written in the
%form:
%$${t_{n-1}}^{k_{n-1}}\cdots {t_1}^{k_1},$$
%and
%$$
%t_i:=s_i s_{i-1}\cdots s_1.
%$$
\end{obs}

This observation was applied in~\cite{AR} to solve a problem of
Foata regarding
 %Problem~\ref{t.foata} for
 the hyperoctahedral group.
%see Propositions~\ref{t.class}(2) and~\ref{B-flag} below.
In this paper, this approach is extended to complex reflection
groups.

\section{Concepts and Results}

\subsection{Background: Wreath Products}

The colored permutation group $G(r,n)$ is the wreath product of
the cyclic group $\bbz_r$ by the symmetric group $S_n$. Namely,
$$
G(r,n) = \bbz_r \wr S_n :=
\{((c_1,\ldots, c_n);\,\pi)\ |\ c_i\in \bbz_r,\,\pi\in S_n\}
$$
with group operation
$$
((c_1,\ldots, c_n);\,\pi) \cdot ((c'_1,\ldots, c'_n);\,\pi') :=
((c_1+c'_{\pi^{-1}(1)},\ldots,c_n+c'_{\pi^{-1}(n)});\,\pi\pi').
$$
%Let
%$$
%\sigma_0:=((1,0,\ldots ,0);\,Id)
%$$
%and
%$$
%\sigma_i:=((0,\ldots,0);\,s_i)\qquad(1\le i< n).
%$$
%The set
%$\{\sigma_0,\sigma_1,\ldots,\sigma_{n-1}\}$ generates $G(r,n)$.

%\vskip 0.2cm

%The following proposition, generalizing Observation~\ref{t.A-fmaj}
%which deals with the special case $G(1,n)=S_n$,
%is a reformulation of the discussion at the beginning of~\cite[Section 2]{AR}.

%\vskip 0.2cm

\begin{pro}\label{t.rn-fmaj}
%Let $\tau_0=Id\in S_n$, $t_i:=s_i s_{i-1}\cdots s_1$ and
%$\tau_i:=((1,0,\dots,0);\ t_i)=t_i s_0=s_i s_{i-1}\cdots s_1s_0$.
%Then
Let $\tau_i := ((1,0,\ldots,0);\,t_i)$ $(0\le i< n)$, where
$t_i := s_i \cdots s_1\in S_n$ $(1\le i< n)$, as in Observation~\ref{t.A-fmaj}
above, and $t_0 = Id\in S_n$, the identity permutation.
Then $(\tau_{n-1},\ldots,\tau_0)$ is a perfect basis for $G(r,n)$, i.e.,
every element $\pi\in G(r,n)$ has a unique presentation
\begin{equation}\label{e.presentation-rn}
\pi={\tau_{n-1}}^{k_{n-1}}\cdots {\tau_1}^{k_1}\tau_0^{k_0},
\end{equation}
where  $0\le k_i < o(\tau_i)=r(i+1)$ $(0\le i< n)$.
\end{pro}

Proposition~\ref{t.rn-fmaj} generalizes the first part of
Observation~\ref{t.A-fmaj}, which concerns the special case $G(1,n)=S_n$.
It is a slightly modified version of a result described in~\cite{AR},
where the basis elements are $\tau_0^{-1} \tau_i \tau_0$ instead of
our $\tau_i$.

Given the unique presentation (\ref{e.presentation-rn}), define the
{\it flag major index} of a colored permutation $\pi\in G(r,n)$ by
$$%\begin{equation}\label{e.fm-rn}
\fmaj_{G(r,n)}(\pi) := \sum_{i=0}^{n-1} k_i,
$$%\end{equation}
the sum of exponents in~(\ref{e.presentation-rn}).

\bigskip

\subsection{General Concepts}

Given a (perfect) basis %(an OGS or an OGS$^*$)
${\bf a}=(a_1,\ldots,a_n)$ for a group $G$, define the
{\it $(G,{\bf a})$ flag major index} as follows.
For every $g\in G$ let
\begin{equation}\label{e.def-fmaj-general}
\fmaj_{(G,{\bf a})}(g) := \sum_{i=1}^n k_i,
\end{equation}
where $k_i$ $(1\le i\le n)$ are the exponents in the unique
presentation
$$
g = a_1^{k_1}\cdots a_n^{k_n} \qquad (0\le k_i< m_i).
$$
Let
$$%\begin{equation}\label{e.fmaj(G)}
\Fmaj_{(G,{\bf a})}(q) := \sum\limits_{g\in G} q^{\fmaj_{(G,{\bf a})}(g)}
$$%\end{equation}
be the corresponding generating function.

By definition,
\begin{equation}\label{e.fmaj(G)-explicit}
\Fmaj_{(G,{\bf a})}(q) = \prod\limits_{i=1}^n [m_i]_q,
\end{equation}
where
$$
[m_i]_q := \frac{q^{m_i}-1}{q-1}.
$$

\bigskip

%\subsection{Mahonian and Hilbertian Order Generating Systems}

Given a group $G$ with a set of generators $S$, let
$\ell_{(G,S)}(\cdot)$ denote the length function on $G$ with
respect to $S$, that is,
$$%\begin{equation}\label{e.ell}
\ell_{(G,S)}(g) :=
\min\{\ell : g = s_1 s_2 \cdots s_\ell \text{ for some }s_i \in S\};
$$%\end{equation}
and let the Poincar\'e series of $G$ (with respect to $S$) be
the corresponding generating function
$$%\begin{equation}\label{e.poincare}
\Poin_{(G,S)}(q) := \sum\limits_{g\in G} q^{\ell_{(G,S)}(g)}.
$$%\end{equation}

The case where $(G,S)$ is a Coxeter system has been extensively
studied  (see, e.g.,~\cite{Hu}). If $G$ is a Coxeter group we will
always assume that $S$ is the Coxeter generating set.
% and abbreviate $\Poin(W,S)$ by $\Poin(W)$.
%This notion generalizes the Poincar\'e series of a Coxeter group
%$W$ (see, e.g.,~\cite[\S 1.11]{Hu}).

\medskip

Motivated by Observation~\ref{t.A-fmaj} we define a
{\em (perfect) Mahonian basis} for $G$ as follows.

\begin{df}\label{mahonian-ogs}
Let ${\bf a}$ be a (perfect) basis for a group $G$ and let $S$ be a %standard
generating set of $G$. Then ${\bf a}$ is a (perfect) {\rm Mahonian basis} for $G$
with respect to $S$ if
$$
\Fmaj_{(G,{\bf a})}(q) = \Poin_{(G,S)}(q);
$$
namely, if the $(G,{\bf a})$ flag major index is equidistributed with length
(with respect to $S$).
\end{df}

\bigskip

Let $V$ be an $n$-dimensional vector space over a field $\bbf$ of
characteristic zero, and let $G$ be a subgroup of the general
linear group $GL(V)$. Then $G$ acts naturally on the symmetric
algebra $S(V^*)$, which may be identified with the polynomial ring
$P_n=\bbf[x_1,\dots,x_n]$. Let $\Lambda^G$ be the subalgebra of
$G$-invariant polynomials, $I^G_n$ the ideal (of $P_n$) generated by
the $G$-invariant polynomials without constant term,
and $R^G := P_n/I^G_n$ the associated coinvariant algebra.
The coinvariant algebra is a direct sum of its homogeneous
components, graded by degree: $R^G = \oplus_k R^G_k$.
%The group $G$ acts naturally on $R^G$ and on each $R^G_k$.
Let
$$%\begin{equation}\label{e.Hilbert}
\Hilb_G(q):= \sum_{k\ge 0} \dim R^G_{k} q^k
$$%\end{equation}
be the corresponding Hilbert series.

%Motivated by Proposition~\ref{t.fm-Hilbert-rn-pro} define a
%{\it (perfect) Hilbertian basis} for $G$ as follows.
\begin{df}\label{d.hilbertian-ogs}
Let ${\bf a}$ be a (perfect) basis for a group $G\subset GL(V)$.
%, whose invariant algebra $\Lambda^G$ is generated by homogenous polynomials.
Then ${\bf a}$ is a (perfect) {\rm Hilbertian basis} for $G$ if
$$%\begin{equation}\label{e.fmaj(G)=Hilb}
\Fmaj_{(G,{\bf a})}(q) = \Hilb_G(q).
$$%\end{equation}
\end{df}

\subsection{Main Result}

Let $r$ be a positive integer and let $p$ be a divisor of $r$. The
complex reflection group $G(r,p,n)$ is defined in~\cite{ST} as the
following subgroup of index $p$ of the wreath product $G(r,n) = \bbz_r \wr S_n$:
$$%\begin{equation}\label{e.grpn}
G(r,p,n) := \{g=((c_1,\ldots,c_n);\ \pi)\in G(r,n)\,|\,
            {\sum_{i=1}^n}c_i\equiv 0 \!\!\pmod p\}.
$$%\end{equation}
For more information on these groups the reader is referred to~\cite{GeckM}.
For the coinvariant algebra and flag major index on these groups see~\cite{BB1}.

The main result of this paper states:

\begin{thm}\label{t.MAIN}
Every complex reflection group $G(r,p,n)$ with parameters satisfying
$gcd(n,p,r/p)=1$ has a perfect Hilbertian basis.
\end{thm}

See Theorem~\ref{t.main} and Corollary~\ref{t.fm-Hilbert-rpn} below.
The special case $p=1$ (wreath product) was established
in~\cite{Stemb, AR1, AR}.

It follows that all classical Weyl groups
have perfect Hilbertian-Mahonian bases (Corollaries~\ref{t.class}
and~\ref{t.BD-flag-equidistribution} below) and that the alternating
subgroup of  a Weyl group of type $B$ has a Mahonian basis
(Proposition~\ref{t.main-AB} below). On the other hand, if
$gcd(n,p,r/p) > 1$ then $G(r,p,n)$ does not necessarily have
a Hilbertian basis; see Section~\ref{s.Appendix} below.

\section{A Perfect Basis for Complex Reflection Groups}

Let ${\bf u}=(u_{n-1},\ldots,u_0)$ be the following sequence of $n$
elements in $G(r,p,n)$:
$$%\begin{equation}\label{e.u-ogs-rpn}
u_i:=(\bc_i;\ t_i)\qquad (0\le i\le n-1),
$$%\end{equation}
where $t_0\in S_n$ is the identity permutation,
$$
t_i := s_i s_{i-1}\cdots s_1 = (i+1,i,\ldots,1)\in S_n\qquad(1\leq i\leq n-1),
$$
$$
\bc_i:=(1,0,\ldots,0,\alpha p-1)\in\bbz_r^n\qquad (0\le i\le n-2),
$$
and
$$
\bc_{n-1}:= (1,0,\ldots,0, p-1).
$$
The integer $0\le \alpha < r/p$ will be chosen later.

\begin{rem}\label{r.beta}
All the results below still hold if we define, more generally,
$$
\bc_{n-1}:= (1,0,\ldots,0, \beta p-1),
$$
where $\beta$ is any integer satisfying $\gcd(\beta,r/p)=1$.
%Note that if $r=1$ then necessarily $\bc_i:=(0,\ldots,0)$ for all $i$.
\end{rem}

\begin{rem}\label{r.zero}
If $r=p$ then one can also take $\bc_{n-1} := (0,\ldots,0)$.
%This will be used later.
\end{rem}

The main result of this section is the following.

\begin{thm}\label{t.main}
If $gcd(n,p,r/p)=1$ then there exists $0\le \alpha< r/p$ such that
${\bf u}$ above is a perfect basis for $G(r,p,n)$.
\end{thm}

%?????
%\vskip 0.2cm
%
%\noindent{\bf Proof.}
%We have to prove that there exists an $\alpha$
%such that every element $\pi\in G(r,p,n)$ has a unique
%presentation of the form:
%
%$$%\begin{equation}\label{ogs-rpn}
%\pi=u_{n-1}^{k_{n-1}}\cdots u_1^{k_1}u_0^{k_0},
%$$%\end{equation}
%where $0\leq k_i< o(u_i) = (i+1)r$ (for $0\leq i\leq n-2$) and
%$0\leq k_{n-1}< o(u_{n-1}) = nr/p$.
%
%!!!!!!

The rest of this section is devoted to proving this result,
using the Chinese Remainder Theorem and the Principle of Inclusion-Exclusion.
For a discussion of the extent to which the condition $\gcd(n,p,r/p)=1$
can be relaxed, see Section~\ref{s.Appendix} below.

\begin{lem}\label{t.isom}
Let $H$ be the subgroup of $G(r,p,n)$ generated by the elements
$\{u_i|\ 0\leq i\leq n-2\}$.
Then $H$ is isomorphic to $G(r,n-1)$.
\end{lem}

\noindent{\bf Proof of Lemma~\ref{t.isom}.}
Define a map $\phi: H \to G(r,n-1)$ by erasing, from each
$\pi = (\bc;\,t)\in H$, the last coordinate of $\bc$.
Let $\psi(\pi)$ be that coordinate, so that $\psi: H \to \bbz_r$.
Since every $\pi\in H$ satisfies $|\pi(n)| = n$, it follows that
$\phi$ and $\psi$ are group homomorphisms. Moreover, for each
$\pi = ((c_1,\ldots,c_{n-1},c_n);\,t)\in H$:
$$
c_n = (\alpha p -1)\sum_{i=1}^{n-1} c_i,
$$
since this property holds for the generators, and is invariant under the
group operation in $H$. It follows that
$$
c_1 = \ldots = c_{n-1} = 0\;\then\;c_n = 0,
$$
namely: $\phi$ is injective. It is also surjective, since
$$
\{((1,0,\ldots,0);\,t_i)\,:\,0\le i\le n-2\}
$$
is a perfect basis for $G(r,n-1)$, by Proposition~\ref{t.rn-fmaj} above.
Thus $\phi$ is a group isomorphism.

\qed

\medskip

Consider now the sequence ${\bf u} = (u_{n-1},\ldots,u_0)$ defined above.
Clearly
$$
o(u_i) = (i+1)r\qquad(0\le i\le n-2)
$$
and
$$
o(u_{n-1}) = nr/p.
$$
(The latter equality holds also if we use the definitions in
Remark~\ref{r.beta} or~\ref{r.zero}.)

The product of all these orders is $n!r^n/p = |G(r,p,n)|$.
If we show that all the products
$$
u_{n-1}^{k_{n-1}} \cdots u_0^{k_0}\qquad(0\le k_i< o(u_i))
$$
are distinct, then it will follow that $\bf u$ is a perfect basis for $G(r,p,n)$.

Assume that
$$
u_{n-1}^{k'_{n-1}} \cdots u_0^{k'_0} = u_{n-1}^{k''_{n-1}} \cdots u_0^{k''_0}
\qquad(0\le k'_i, k''_i < o(u_i)).
$$
We want to show that $k'_i = k''_i$ $(\forall i)$.
It suffices to show that $k'_{n-1} = k''_{n-1}$, since then
$$
u_{n-2}^{k'_{n-2}} \cdots u_0^{k'_0} = u_{n-2}^{k''_{n-2}} \cdots u_0^{k''_0}
$$
and, by Lemma~\ref{t.isom} and Proposition~\ref{t.rn-fmaj},
this implies $k'_i = k''_i$ $(0\le i\le n-2)$.

Indeed, by assumption
$$
u_{n-1}^{k'_{n-1}-k''_{n-1}} =
[u_{n-2}^{k''_{n-2}}\cdots u_0^{k''_0}]
[u_{n-2}^{k'_{n-2}}\cdots u_0^{k'_0}]^{-1}\in H.
$$
Let $k := k'_{n-1}-k''_{n-1}$; working modulo $o(u_{n-1})$,
we can assume that $0\le k < nr/p$.
$u_{n-1}^k\in H$ implies that $|u_{n-1}^k(n)| = n$ and therefore,
by considering the $S_n$-component of $u_{n-1}$, $n|k$.
Denote $\tk := k/n$. Then
$$%\begin{equation}\label{e.uk}
u_{n-1}^k = u_{n-1}^{n\tk} = ((\tk p,\ldots,\tk p);Id),
$$%\end{equation}
where $Id\in S_n$ is the identity permutation and $0\le \tk < r/p$.
(If we use the definition in Remark~\ref{r.beta} then $\tk p$ should be
replaced here by $\tk \beta p$. If we use the definition in
Remark~\ref{r.zero} then $o(u_{n-1}) = nr/p = n$, and the proof here.)

On the other hand, we can present $u_{n-1}^k\in H$ in the form
$$
u_{n-1}^k = u_{n-2}^{k_{n-2}}\cdots u_0^{k_0}\qquad(0\le k_i < o(u_i)).
$$
The natural projection $T: H \to S_{n-1}$, defined by $T((\bc;\,t)) := t$,
is a group homomorphism mapping the perfect basis $(u_{n-2},\ldots,u_0)$
of $H$ onto the perfect basis $(t_{n-2},\ldots,t_0)$ of $S_{n-1}$.
Since $T(u_{n-1}^k) = Id$, it follows that $o(t_i) = i+1$ divides $k_i$;
let $\tk_i := k_i/(i+1)$ ($0\le i\le n-2$).
Now
$$
u_i^{i+1} = (v_i;Id)\qquad(0\le i\le n-2)
$$
where
\begin{equation}\label{e.ui}
v_i := (\underbrace{1,\ldots,1}_{i+1},0,\ldots,0,(\alpha p - 1)(i+1))\in\bbz_r^n.
%\qquad(0\le i\le n-2).
\end{equation}
Thus
$$
u_{n-1}^k = u_{n-2}^{(n-1)\tk_{n-2}}\cdots u_1^{2\tk_1}u_0^{\tk_0}
= (\sum_{i=0}^{n-2}\tk_i v_i;Id).
$$
So far we have
$$
\sum_{i=0}^{n-2}\tk_i v_i = (\tk p,\ldots,\tk p) \in \bbz_r^n\
\qquad(0\le \tk_i < \frac{o(u_i)}{i+1}=r,\,0\le i\le n-2).
$$
Since $v_0,\ldots,v_{n-2}\in\bbz_r^n$ are linearly independent, we conclude that
$$
\tk_{n-2} = \tk p
$$
while
$$
\tk_i = 0\qquad(0\le i\le n-3).
$$
Thus
$$
\tk p v_{n-2} = (\tk p, \ldots, \tk p).
$$
Comparing the last coordinate on each side, we get by%~(\ref{e.uk}) and
~(\ref{e.ui}):
$$
\tk p (\alpha p - 1)(n-1) = \tk p\qquad(\hbox{in\ }\bbz_r).
$$
(Multiply both sides by $\beta$ for Remark~\ref{r.beta}.)
Rewriting $(\alpha p - 1)(n-1) - 1 = (n-1)\alpha p - n$, this is equivalent to
\begin{equation}\label{e.mp}
%\tk[(\alpha p - 1)(n-1) - 1] = 0\qquad(\hbox{in\ }\bbz_{r/p}),
\tk [(n-1)\alpha p - n] = 0\qquad(\hbox{in\ }\bbz_{r/p}),
\end{equation}
where $0\le \tk < r/p$ and $0\le \alpha < r/p$.
(Same equation for Remark~\ref{r.beta}, since $\gcd(\beta,r/p)=1$.)
%We want to prove: $\tk=0$.

We want to show that there exists $0\le \alpha < r/p$ such that
(\ref{e.mp}) necessarily implies $\tk = 0$.
Equivalently, we must find $\alpha$ such that
%$r/p$ does not divide
%$\tk[(\alpha p - 1)(n-1) - 1] = \tk[(n-1)\alpha p - n]$,
%for any $1\le \tk < r/p$. In other words:
$$%\begin{equation}\label{e.gcd}
\gcd(r/p,(n-1)\alpha p - n)=1.
$$%\end{equation}
If $r/p = 1$, every $\alpha$ will do. In general, we want to show that
the following ``False Assumption'' leads to a contradiction.

\vskip 0.2cm

\noindent
%\begin{itemize}\item[(*)]
{\bf False Assumption:} For every $0\le \alpha < r/p$,
$$
\gcd(r/p, (n-1)\alpha p-n) > 1.
$$
%there exists a divisor $q_\alpha > 1$ of $r/p$ which
%divides $(n-1)\alpha p-n$.
%\end{itemize}

\vskip 0.2cm

%The case $\frac{r}{p}=1$ holds trivially, since there are no
%divisors greater than 1 of $\frac{r}{p}=1$ in this case. For the
%general case we need the following lemma.

%\vskip 0.2cm

\begin{lem}\label{t.1}%\rm
If $q>1$ is a common divisor of $r/p$, $(n-1){\alpha}p-n$ and
$(n-1){\beta}p-n$, where $\alpha \ne \beta$ and $gcd(\beta-\alpha,q)=1$,
then $q$ divides $\gcd(n,p,r/p)$.
\end{lem}

%\vskip 0.2cm

\noindent {\bf Proof of Lemma~\ref{t.1}.}
By assumption, $q$ divides
$((n-1){\beta}p-n) - ((n-1){\alpha}p-n) = (n-1){(\beta-\alpha)}p$.
Since $gcd(\beta-\alpha,q)=1$, $q$ divides $(n-1)p$. Thus $q$ divides
$\alpha(n-1)p - ((n-1){\alpha}p-n) = n$, so that $gcd (q,n-1)=1$.
Hence $q$ divides $p$ as well, completing the proof of the lemma.

\qed

\vskip 0.2cm

By the ``False Assumption'' above there exists, for each $0\le \alpha < r/p$,
a common (prime) divisor of $r/p$ and $(n-1)\alpha p-n$.

%\vskip 0.2cm

\begin{lem}\label{t.q}%\rm
Assume that $\gcd(n, p, r/p) = 1$, and denote
$$
Q := \{q\hbox{ prime}\,|\, q\hbox{ divides } r/p \hbox{ and } (n-1){\alpha}p-n
                           \hbox{ for some } 0\le \alpha < r/p\}.
$$
Then, for any number of distinct primes $q_1, \ldots, q_t\in Q$, the number of
integers $0\le \alpha < r/p$ such that $(n-1)\alpha p-n$ is divisible by all
of $q_1, \ldots, q_t$ is $r/(p q_1 \cdots q_t)$.
\end{lem}

%\vskip 0.2cm

\noindent {\bf Proof of Lemma~\ref{t.q}.}
Let $q\in Q$, and assume that it divides $(n-1){\alpha}p-n$.
If $\beta - \alpha$ is divisible by $q$, then clearly $q$ divides also
$(n-1){\beta}p-n$. Conversely, if $\beta - \alpha$ is not divisible by the
prime $q$ then $\gcd(\beta - \alpha,q) = 1$. By Lemma~\ref{t.1}, and since
$\gcd(n,p,r/p) = 1$ by assumption, $q$ does not divide $(n-1){\beta}p-n$.
It follows that the number of $0\le \alpha < r/p$ divisible by any $q\in Q$
is exactly $r/(pq)$.

We now consider any number of distinct primes
$q_1,\ldots,q_t\in Q$. Suppose that
$q_i$ divides $(n-1){\alpha_i}p-n$ $(1\le i\le t)$.
By the above argument, an integer $\alpha$ has the property that
$(n-1){\alpha}p-n$ is divisible by all of the $q_i$ if and only if
$\alpha$ solves the $t$ simultaneous modular equations
$$%\begin{equation}\label{e.mod}
\alpha \equiv \alpha_i \pmod{q_i}\qquad(1\le i\le t).
$$%\end{equation}
A solution exists, and is unique $\!\!\pmod{q_1 \cdots q_t}$,
by the Chinese Remainder Theorem.
It follows that the number of $0\le \alpha < r/p$ divisible by all of
$q_1, \ldots, q_t$ is exactly $r/(p q_1 \cdots q_t)$.

\qed

\vskip 0.2cm

We shall now wrap up, by counting the integers $0\le \alpha < r/p$ according
to which primes $q\in Q$ divide $(n-1){\alpha}p-n$.
According to the ``False Assumption'', each $\alpha$ has at least one such $q$.
By Lemma~\ref{t.q} and the Principle of Inclusion-Exclusion, counting gives
$$
\frac{r}{p} = \sum_{q\in Q} \frac{r}{pq} - \sum_{q_1 < q_2} \frac{r}{p q_1 q_2}
+ \sum_{q_1 < q_2 < q_3} \frac{r}{p q_1 q_2 q_3} - \ldots.
$$
Rearrangement gives
$$
\frac{r}{p} \cdot \prod_{q\in Q} \left(1 - \frac{1}{q}\right) = 0,
$$
which is clearly a contradiction, since $Q$ is a finite set
of integers greater than $1$.
This completes the proof of Theorem~\ref{t.main}.

\section{Identities}

\subsection{A Flag Major Index %and a Perfect Hilbertian Basis
for $G(r,p,n)$}

$G(r,p,n)$ is a subgroup of $G(r,n)$ and thus acts naturally on
the polynomial ring $P_n=\bbq[x_1,\dots,x_n]$. Denote the ring of
$G(r,p,n)$-invariant polynomials in $P_n$ by $\Lambda_{r,p,n}$.
Let $I_{r,p,n}$ be the ideal of $P_n$ generated by the elements
of $\Lambda_{r,p,n}$ without constant term. The quotient
$R_{r,p,n} := P_n/I_{r,p,n}$ is the {\it coinvariant algebra} of
$G(r,p,n)$. Each complex reflection group $G(r,p,n)$ acts
naturally on its coinvariant algebra. Let $R^{(k)}_{r,p,n}$ be the
$k$-th homogeneous component of the coinvariant algebra,
$R_{r,p,n}=\oplus_k R^{(k)}_{r,p,n}$, and let
$$
\Hilb_{r,p,n}(q) := \sum_{k\ge 0} \dim R^{(k)}_{r,p,n} q^k
$$
be the corresponding Hilbert series.
%The Hilbert series of the coinvariant algebra of $G(r,p,n)$,
$\Hilb_{r,p,n}(q)$ was expressed in~\cite{BB1} as a generating function
for $\fmaj_{G(r,n)}$ on a certain subset of the wreath product $G(r,n)$.
Using Theorem~\ref{t.main}
%we define a natural flag major index on $G(r,p,n)$ and show that
it will be shown that $\Hilb_{r,p,n}(q)$ may be expressed as a
generating function for a natural flag major index on the group
$G(r,p,n)$ itself. This generalizes results for $G(r,1,n)$ which
were proved in~\cite{Stemb, AR1, AR}.
%identities (\ref{e.fm-Hilbert-S_n}) and
%(\ref{e.fm-Hilbert-rn}).

\smallskip

Let $G := G(r,p,n)$ with $gcd(n,p,r/p)=1$.  Recall the
perfect basis {\bf u} for $G$ from Theorem~\ref{t.main} and the flag
major index $\fmaj_{(G, {\bf u})}$  from
Definition~(\ref{e.def-fmaj-general}).
%Then for $\pi\in G(r,p,n)$
%\begin{equation}\label{def-fmaj-rpn}
%\fmaj_{(G(r,p,n), {\bf u})}(\pi):=\sum\limits_{i=0}^{n-1} k_i,
%\end{equation}
%where $\pi=u_{n-1}^{k_{n-1}}\cdots u_1^{k_1}u_0^{k_0}$, $u_i$ are
%defined as in (\ref{e.u-ogs-rpn}) and $\alpha$ is chosen as in the
%proof of Theorem~\ref{t.main}.

%The following corollary generalizes  identities
%(\ref{e.fm-Hilbert-S_n}) and (\ref{e.fm-Hilbert-rn}).

\begin{cor}\label{t.fm-Hilbert-rpn}
If $gcd(n,p,r/p)=1$ then {\bf u} is a perfect Hilbertian
basis for $G(r,p,n)$; namely,
$$
\Hilb_{r,p,n}(q) = \Fmaj_{(G(r,p,n), {\bf u})}(q),
$$
where $\Fmaj_{(G(r,p,n), {\bf u})}(q) :=
\sum\limits_{\pi \in S_n} q^{\fmaj_{(G(r,p,n), {\bf u})}(\pi)}$.
\end{cor}

%This result should be compared to~\cite{BB1}, where $\Hilb
%(r,p,n)$ was expressed as a generating function of
%$\fmaj_{G(r,n)}$ on a subset of the wreath product $G(r,n)$.

\noindent{\bf Proof.}
By Theorem~\ref{t.main} and identity (\ref{e.fmaj(G)-explicit}),
$$
\sum\limits_{\pi \in S_n} q^{\fmaj_{G(r,p,n)}(\pi)} =
[r]_q [2r]_q \cdots [(n-1)r]_q [nr/p]_q
$$
where $[m]_q := \frac{q^m-1}{q-1}$.
On the other hand, it is well known (see, e.g.,~\cite{BB1}) that
\begin{equation}\label{e.Hilb-comp}
\Hilb_{r,p,n}(q) = [r]_q [2r]_q \cdots [(n-1)r]_q [nr/p]_q,
\end{equation}
completing the proof.
%Equation (\ref{e.Hilb-comp}) completes the proof.
\qed

\subsection{Classical Weyl Groups}

Recall the three infinite series of classical Weyl group: the
symmetric groups $S_n$ (Weyl groups of type $A$), the signed
permutation groups (sometimes called hyperoctahedral groups) $B_n$
(Weyl groups of type $B$), and the even signed permutation groups
$D_n$ (Weyl groups of type $D$).

\begin{cor}\label{t.class}$\,$
\begin{itemize}
\item[1.]
Let
$$
\alpha_i := (i,i-1,\ldots,1) = [i,1,2,\ldots,i-1,i+1,i+2,\ldots,n]
\qquad(2\le i\le n)
$$
be permutations.
The sequence ${\bf a}=(\alpha_n, \alpha_{n-1},\ldots,\alpha_2)$
is a perfect Hilbertian basis %ordered generating system
for the symmetric group $S_n$.
\item[2.]
Let
$$
\beta_i := [-i,1,2,\ldots,i-1,i+1,i+2,\ldots,n]
\qquad(1\le i\le n)
$$
be signed permutations.
The sequence ${\bf b}=(\beta_n, \beta_{n-1},\ldots,\beta_1)$
is a perfect Hilbertian basis % ordered generating system
for the hyperoctahedral group $B_n$.
\item[3.]
Let
$$
\delta_i:=[-i,1,2,\ldots,i-1,i+1,i+2,\ldots,-n]\qquad(1\le i\le n-1)
$$
and
$$
\delta_n:=[n,1,2,\ldots,n-1]
$$
be signed permutations.
The sequence ${\bf d}=(\delta_n, \delta_{n-1},\ldots,\delta_1)$
is a perfect Hilbertian basis %ordered generating system
for the group of even signed permutations $D_n$.
\end{itemize}
\end{cor}

\noindent{\bf Proof.}
By Theorem~\ref{t.main} and Remark~\ref{r.zero}, ${\bf a}$, ${\bf b}$
and ${\bf d}$ are perfect bases for $S_n=G(1,1,n)$, $B_n=G(2,1,n)$
and $D_n=G(2,2,n)$, respectively
(using $\alpha = 1$ for $B_n$ and $D_n$, with Remark~\ref{r.zero} for $D_n$).
By Corollary~\ref{t.fm-Hilbert-rpn}, these bases are Hilbertian. \qed

\begin{cor}\label{t.BD-flag-equidistribution}
\begin{itemize}
\item[1.] The sequence ${\bf b}$ is a perfect Mahonian basis for
$B_n$ (with respect to the Coxeter generating set  $S$). Namely,
the resulting flag major index $\fmaj_{(B_n, {\bf b})}$ is
equidistributed with the length function $\ell_{(B_n,S)}$ over
$D_n$.
% For $\pi\in B_n$,
%let
%$$\fmaj_{(B_n, {\bf b})}(\pi):=\sum_{i=1}^{n}k_i,$$ where $\pi$ is written in the form:
%${\beta_n}^{k_n}\cdots {\beta_1}^{k_1}$, and  $0\le i<2i$ for
%$1\le i\le n$.
%Then
%\begin{equation}\label{B-equidistribution}
%\sum_{\pi \in B_n}q^{\ell_{B_n}(\pi)} = \sum_{\pi \in B_n}
%q^{\fmaj_{B_n}(\pi)},
%\end{equation}
%where $\ell_{B_n}$ is the standard length function with respect to
%the Coxeter generating set of $B_n$.

\item[2.] The sequence ${\bf d}$ is a perfect Mahonian basis for
$D_n$ (with respect to the Coxeter generating set $S'$). Namely,
the resulting flag major index $\fmaj_{(D_n, {\bf d})}$ is
equidistributed with the length function $\ell_{(D_n,S')}$ over
$D_n$.
%For $\pi\in
%D_n$, let
%$$\fmaj_{(D_n, {\bf d})}(\pi):=\sum_{i=1}^{n}k_i,$$ where $\pi$ is written in the form:
%$$
%{{\delta_n}^{k_n}\delta_{n-1}}^{k_{n-1}}\cdots
%{\delta_1}^{k_1}\qquad 0\le k_i \le 2i \hbox{ for } 1\le i< n
%\hbox{ and } 0\le k_n <n .
%$$
%Then
%\begin{equation}\label{D-equidistribution}
%\sum_{\pi \in D_n}q^{\ell_{(D_n, S')}(\pi)} = \sum_{\pi \in D_n}
%q^{\fmaj_{D_n}(\pi)},
%\end{equation}
%where $\ell_{(D_n, S')}$ is the standard length function with respect to
%the Coxeter generating set of $D_n$.

\end{itemize}
\end{cor}

\noindent{\bf Proof.} It is well known that for every Weyl group
$W$, the Hilbert series of the coinvariant algebra of $W$ is equal
to the Poincar\'e series of $W$, namely to the generating function
for length with respect to the Coxeter generators; see, e.g.,
\cite[\S 3.15]{Hu}. Combining this with Corollary~\ref{t.fm-Hilbert-rpn}
gives the desired result. \qed

\medskip

%\begin{rem}
While the statements on types $A$ and $B$
%(Corollary~\ref{t.class}(1)-(2) and
%equation~(\ref{B-equidistribution}))
are not new, see~\cite{AR}, the statements on type $D$
(Corollary~\ref{t.class}(3) and
Corollary~\ref{t.BD-flag-equidistribution}(2))
%equation~(\ref{D-equidistribution}))
are new. In particular, note that $\fmaj_{(D_n, {\bf u})}$ is
equidistributed with, but different from, the flag major index for $D_n$
which was introduced by Biagioli and Caselli~\cite{BC1}.
%\end{rem}

\subsection{The Alternating Subgroup of $B_n$}

%\vskip 0.2cm

%{\bf Definition:}

Let $B_n^+$ be the alternating subgroup of the Coxeter group of type
$B$; namely,  the subgroup consisting all elements in $B_n$  of
even length.

Let $r_1:=[2,-1,3,\dots,n]$ and
$r_i:=[-1,2,\dots , i+1,i, i+2, i+3, \dots, n]$ ($2\le i\le n-1$).
$R:=\{r_i\ |\ 1\le i\le n-1\}$ is a set of generators for $B_n^+$ with
Coxeter-like relations~\cite[Chapter IV Section 1 Exercise 9]{Bou}.
The defining relations are:
$$
r_1^4=1
$$
$$
r_i^2=1 \qquad (1<i< n)
$$
$$
(r_ir_{i+1})^3=1 \qquad (1\le i< n).
$$
$$
(r_ir_j)^2=1 \qquad (|i-j|>1)
$$
Let $\ell_{(B_n^+, R\cup R^{-1})}(\pi)$ be the length of $\pi\in B_n^+$
with respect to $R\cup R^{-1}$. Let

%Then
%
%\begin{pro}\label{t.length-ABn}~\cite[Example 5.2.6]{BRR}
%$$\sum\limits_{\pi \in B_n^+} q^{\ell(\pi)}=\prod_{i=1}^{n-1}[2i]_q[n]_q.$$
%%where $[m]_q:={q^m-1\over q-1}$.
%\end{pro}
%
%%\vskip 0.2cm

$$
v_n := %s_{n-1} \cdots s_1s_0s_1 \cdots s_{n-1}=
((0,\ldots,0,1);\ Id)
     = [1,2,\ldots,-n]\in B_n,
$$
and define a map $\psi:D_n\mapsto B_n^+$ by
$$
\psi(w):=   \begin{cases}
      w    & \text{ if } w \in B_n^+ ; \\
      wv_n & \text{ if } w \not\in B_n^+ .
   \end{cases}
$$
Namely, $\psi$ switches the sign of the last letter of $w$ if
$w\not\in B_n^+$.

\begin{fac}\label{bijection1}
$\psi$ is a bijection.
\end{fac}

Recall the basis ${\bf d}=(\delta_n,\dots,\delta_1)$ for $D_n$
from Corollary~\ref{t.class}(3) and let
$$
\gamma_i:=\psi(\delta_i)\qquad (1\le i\le n).
$$

\begin{pro}\label{t.main-AB}
%The sequence $(\gamma_n,\dots,\gamma_1)$ is an ordered generating
%system for $B_n^+$.
(1). The sequence ${\bf c}=(\gamma_1,\dots,\gamma_n)$ is a
Mahonian basis for $B_n^+$. Namely
\begin{itemize}
\item[(i)]
Every element $\pi\in B_n^+$ has a unique presentation
\begin{equation}\label{AB-presentation}
 \pi=\gamma_n^{k_n}\gamma_{n-1}^{k_{n-1}}\cdots
\gamma_1^{k_1}\qquad  0\le k_i \le 2i \hbox{ for } 1\le i< n
\hbox{ and } 0\le k_n <n .
\end{equation}
\item[(ii)]
\begin{equation}\label{AB-equidistribution}
\sum\limits_{\pi\in B_n^+} q^{\fmaj_{(B_n^+, {\bf
c})}(\pi)}=\sum\limits_{\pi\in B_n^+} q^{\ell_{(B_n^+, R\cup
R^{-1})}(\pi)}.
\end{equation}
\end{itemize}
\noindent
(2).  % For every $\pi\in B_n^+$ let
%$$
%\fmaj_{(B_n^+, {\bf c})}(\pi):=\sum\limits_{i=1}^n k_i,
%$$
%where $k_1,\dots,k_n$ are the exponents in the presentation
%(\ref{AB-presentation}). Then t
The flag major index is invariant under $\psi$. Namely, for every
$w\in D_n$
\begin{equation}\label{e.flag-invariance}
\fmaj_{(D_n, {\bf d})}(w)=\fmaj_{(B_n^+, {\bf c})} (\psi(w)).
\end{equation}
\noindent
(3). For every $w\in D_n$,
$\fmaj_{(D_n, {\bf d})}(w)\equiv {0 \pmod 2}$
if and only if $w\in D_n\cap B_n^+$.
Similarly, for every $w\in B_n^+$,
$\fmaj_{(B_n^+, {\bf c})}(w)\equiv {0 \pmod 2}$
if and only if $w\in D_n\cap B_n^+$.
\end{pro}

\noindent{\bf Proof.}
Let $w$ be an element in $D_n$. By
Corollary~\ref{t.class}(3), there exist unique $0\le k_i <2i$ $(0\le
i< n$) and $0\le k_n<n$ such that $w=\delta_n^{k_n}\cdots
\delta_1^{k_1}$. Noticing that $v_n$ commutes with $\delta_i$ for
$i<n$ we obtain
$$
\gamma_n^{k_n}\cdots \gamma_1^{k_1}=(\delta_n v_n)^{k_n}\cdots
(\delta_1 v_n)^{k_1}=\delta_n^{k_n}\cdots \delta_1^{k_1}
v_n^{\sum_i k_i} $$ $$= w v_n^{\fmaj_{(D_n, {\bf d})}(w)}=w
v_n^{\fmaj_{(D_n, {\bf d})}(w)\ {\rm mod}\ 2}.
$$
But $\gamma_n^{k_n}\cdots \gamma_1^{k_1}\in B_n^+$ while
$v_n\not\in B_n^+$. It follows that $w\in B_n^+$ if and only if
$\fmaj_{(D_n, {\bf d})}(w)\ {\rm mod}\ 2= 0$. Hence
$$
\gamma_n^{k_n}\cdots \gamma_1^{k_1}=w v_n^{\fmaj_{(D_n, {\bf
d})}(w)\ {\rm mod}\ 2}=\psi(w).
$$
Since $\psi$ is a bijection this proves $(i)$, (2) and (3).

\smallskip

To prove $(ii)$ recall from~\cite{BRR} the bijection
$\theta:B_n^+\mapsto D_n$
$$
\theta(w):=\begin{cases}
      w    & \text{ if } w \in B_n^+ ; \\
      w s_0 & \text{ if } w \not\in B_n^+,
   \end{cases}
$$
which switches the sign of the first letter of $w$ if $w\not\in
D_n$. By~\cite[Corollary 5.2(i)]{BRR}, the length is invariant
under $\theta$. Namely, for every $w\in B_n^+$
\begin{equation}\label{e.length-invariance}
\ell_{(B_n^+, R\cup R^{-1})}(w)=\ell_{(D_n, S')}(\theta(w)).
\end{equation}
Combining (\ref{e.flag-invariance}), (\ref{e.length-invariance}) with
Corollary~\ref{t.BD-flag-equidistribution}(2) and the fact that
$\psi$ and $\theta$ are bijections we obtain
$$
\sum\limits_{\pi\in B_n^+} q^{\fmaj_{(B_n^+, {\bf c})}(\pi)} =
\sum\limits_{\pi\in B_n^+} q^{\fmaj_{(D_n, {\bf d})}(\psi^{-1}(\pi))} =
\sum\limits_{w\in D_n} q^{\fmaj_{(D_n, {\bf d})}(w)} =
$$
$$
=
\sum\limits_{w\in D_n} q^{\ell_{(D_n, S')}(w)} =
\sum\limits_{w\in D_n} q^{\ell_{(B_n^+, R\cup R^{-1})}(\theta^{-1}(w))} =
\sum\limits_{\pi\in B_n^+} q^{\ell_{(B_n^+, R\cup R^{-1})}(\pi)}.
$$
This completes the proof of $(ii)$.

\qed

\bigskip

\noindent{\bf Remarks.}
1. $(\gamma_n,\dots,\gamma_1)$ is a
perfect Mahonian basis for $B_n^+$ if and only if $n$ is odd.
If $n$ is even then the order of $\gamma_n$ is $2n$ while $k_n$ is bounded
by $n$ in~(\ref{AB-presentation}), so $B_n^+$ is not decomposed
into a set-wise direct product of the cyclic subgroups generated
by $\gamma_n,\ldots,\gamma_1$; in this case $(\gamma_n,\ldots,\gamma_1)$ is a
Mahonian basis for $B_n^+$ which is not perfect. %An ordered generating set which gives a
%unique expression of type (\ref{AB-presentation}) but does not
%give a decomposition into a set-wise direct product of the cyclic
%subgroups is called an OGS$^*$, see~\cite{S1}.

\smallskip

\noindent
2. A major index and a Mahonian identity on the alternating subgroup of $S_n$
may be found in~\cite{RegevR}. It should be noted that, while the length
function is defined there with respect to a generating set analogous to
the above $R\cup R^{-1}$, there is apparently no simple interpretation,
involving bases, of the major index in this case.

\section{Complex Reflection Groups with No Hilbertian Basis}\label{s.Appendix}

\begin{pro}\label{t.no_hilbert}
For any prime $p$, the group $G(p^2,p,p)$ has no perfect Hilbertian basis.
\end{pro}

%\vskip 0.2cm

\noindent{\bf Proof.}
Assume that $p$ is a prime number for which $G(p^2,p,p)$ has a perfect
Hilbertian basis. A Hilbert function of the form~(\ref{e.Hilb-comp}) has
a unique decomposition into factors of the form $[m_i]_q$, where $m_i$
are positive integers. It follows that, up to reordering, the $p$ elements
$t_0, t_1, \ldots, t_{p-1}$ in a perfect Hilbertian basis for $G(p^2,p,p)$
have orders
$o(t_0)=p^2$, $o(t_1)=2p^2$, $\ldots, o(t_{p-2})=(p-1)p^2$ and
$o(t_{p-1})=p^2$.
Let $t_i = (v_i;\ \pi_i)$, where $v_i\in (\bbz_{p^2})^p$ with sum of entries
$\equiv 0 \pmod p$ and $\pi_i\in S_p$ $(0\le i\le p-1)$.

Both $t_0$ and $t_{p-1}$ are of order $p^2$, and therefore
neither $\pi_0$ nor $\pi_{p-1}$ contains a cycle of any size $1<i<p$.
Each of them is, therefore, either a $p$-cycle or the identity permutation.
If ${\pi}_0$ is a $p$-cycle then $t_0^p=(w_0;\ Id)$ where
$w_0=(\alpha,\ldots,\alpha)$,
$\alpha\equiv 0 \pmod p$ but $\alpha\not\equiv 0 \pmod {p^2}$.
If ${\pi}_0=Id$ then $t_0^p=(w_0;\ Id)$ where all the entries of $w_0$ are
$0 \pmod p$ but not all are $0 \pmod {p^2}$, and their sum is $0 \pmod {p^2}$.
In both cases, $w_0\in (p\bbz_{p^2})^p$ is a nonzero vector whose sum of
entries is $0 \pmod {p^2}$. The same conclusion holds for $w_{p-1}$, where
$t_{p-1}^p=(w_{p-1};\ Id)$

Now let $1 < i < p$.
Then $o(t_{i-1})=ip^2$, and therefore $o(\pi_{i-1})\,|\,ip^2$.
If $\pi_{i-1}$ is a $p$-cycle then $t_{i-1}^{p^3}=Id$; and since
$\gcd(i,p) = 1$ this implies $t_{i-1}^{p^2}=Id$, contradicting $1 < i < p$.
Thus $\pi_{i-1}$ is not a $p$-cycle, and therefore $o(\pi_{i-1})\,|\,i$.
Denoting $t_{i-1}^{ip} = (w_{i-1};\ Id)$ $(1 < i < p)$, it follows that
$w_{i-1}\in (p\bbz_{p^2})^p$ is
a nonzero vector whose sum of entries is $0 \pmod {p^2}$.

We conclude that all the vectors $w_0, w_1, \ldots, w_{p-2}, w_{p-1}$
belong to
$$
V := \{w = (\alpha_1,\ldots,\alpha_p)\in (p\bbz_{p^2})^p\,|\,
\alpha_1 + \ldots + \alpha_p = 0\},
$$
which is a $(p-1)$-dimensional vector space over the field $\bbz_p$.
The unique presentation property of the basis $t_0, \ldots, t_{p-1}$
implies that these $p$ vectors are linearly independent over $\bbz_p$.
This is a contradiction which completes the proof of the proposition.
\qed

\medskip

\end{document}